\input amstex
\documentstyle{amsppt}
\topmatter
\title
A Decomposition Theorem on Differential Polynomials of
Theta Functions of High Level
\endtitle
\rightheadtext{Theta Functions of High Level}
\author  Jae-Hyun Yang
\endauthor
\magnification=\magstep 1 \baselineskip =7mm \pagewidth{12.5cm}
\pageheight{18.00cm}
\thanks{This work was in part supported by TGRC-KOSEF}
\endthanks
\address{Department of Mathematics \endgraf Inha University\endgraf
Incheon 402-751 \endgraf Republic of Korea\endgraf
email\,:\,jhyang\@inha.ac.kr}
\endaddress
\endtopmatter
\document
\NoBlackBoxes

\define\pw{\left({{\partial}\over {\partial W}}\right)}
\define\bhg{\BZ^{(h,g)}_{\geq 0}}
\define\BM{{\Bbb M}(h)}
\define\LM{{\Cal L}_{\Cal M}}
\define\lrt{\longrightarrow}
\define\lmt{\longmapsto}

\define\s{\sigma}
\define\BZ{\Bbb Z}

\define\BR{\Bbb R}
\define\BC{\Bbb C}
\define\J{J\in {\Bbb Z}^{(h,g)}_{\geq 0}}
\define\N{N\in {\Bbb Z}^{(h,g)}}

\define\M{\Cal M}
\define\ta{\vartheta^{(\Cal M)}\left[\matrix A\\ 0\endmatrix\right]
(\Omega\vert W)}
\define\Jt{\left({{\partial}\over {\partial W}}\right)^J\vartheta^{(\Cal M)}
\left[\matrix A\\ 0\endmatrix\right](\Omega\vert W)}
\define\Dt{\Delta^J\vartheta^{(\Cal M)}\left[\matrix A\\ 0\endmatrix\right]
(\Omega\vert W)}
\define\emw{\text {exp}\left\{-\pi i\sigma({\Cal M}(\xi\Omega\,^t\xi
+2W\,^t\!\xi))\right\} }
\define\mjz{{\tilde {\vartheta}}^{(\Cal M)}_J\left[\matrix A\\ 0\endmatrix
\right](\Omega\vert Z,W)}
\define\CJ{{\Bbb C}\left[\cdots,\left({{\partial}\over
{\partial W}}\right)^J \vartheta^{(\Cal M)}\left[\matrix
A_{\alpha,{\Cal M}}\\ 0\endmatrix \right](\Omega\vert
W),\cdots\right]}
\define\GJ{G\left(\cdots,\Delta^J\vartheta^{(\Cal M)}\left[\matrix
A_{\alpha,{\Cal M}}\\ 0\endmatrix\right](\Omega\vert
W),\cdots\right)}
\define\GW{G\left(\cdots,\left({{\partial}\over {\partial W}}\right)^J
\vartheta^{(\Cal M)}\left[\matrix A_{\alpha,{\Cal M}}\\
0\endmatrix\right] (\Omega\vert W),\cdots\right)}

\ \ \  Let $h$ and $g$ be two positive integers. We fix an element
$\Omega$ of the Siegel upper half plane
$$H_g:=\left\{\,Z\in \BC^{(g,g)}\,\vert\ Z=\,^tZ,\ \ \text {Im}\,Z>0\
\right\}$$
of degree $g$ once and for all. Let $\M$ be positive symmetric,
even integral matrix of degree $h$. An entire function $f$ on $\BC^{(h,g)}$
satisfying the transformation behaviour
$$F(W+\xi \Omega+\eta)=\emw\,f(W)$$
for all $W\in \BC^{(h,g)}$ and $(\xi,\eta)\in \BZ^{(h,g)}\times
\BZ^{(h,g)}$ is called a {\it theta\ function\ of\ level} $\M$ with respect
to $\Omega.$ The set $T_{\M}(\Omega)$ of all theta functions of level $\M$
with respect to $\Omega$ forms a complex vector space of dimension
$\left(\,\text {det}\,\M\,\right)^g$ with a canonical basis consisting of
theta series
$$\ta:=\sum_{\N}\,\text {exp}\left\{\pi i\s(\M((N+A)\Omega\,^t(N+A)+
2W\,^t(N+A))))\right\},$$
where $A$ runs over a complete system of representatives of the cosets
$\LM:=\M^{-1}\BZ^{(h,g)}/\BZ^{(h,g)}.$ \par
\ \ \ We let
$$T(\Omega):=\sum_{\M}\,T_{\M}(\Omega)$$
be the graded algebra of theta functions, where $\M=(\M_{kl})\,
(1\leq k,l\leq h)$ runs over the set $\BM$ of all positive symmetric,
even integral $h\times h$ matrices with $\M_{kl}\not= 0$ for all $k,l.$
\vskip 0.2cm
\ \ \
In this paper we prove the following decomposition theorem:\par
The algebra of differential polynomials of theta functions has a canonical
basis
$$\left\{\,\Jt\,\bigg|\ \J,\ A\in \LM,\ \M\in \BM\,\right\}, $$
i.e., any differential polynomials of theta functions can be expressed
uniquely as a linear combination of $\Jt\ (\,\J,\,A\in \LM,\,\M\in \BM\,)$
with constant coefficients depending only on $\Omega.$ \par
\ \ \ The key idea is a quite
similiar one as making transvectants in the classical invariant theory\,(cf.\,
[M1],[M2]\,). However the Lie algebra is the Heisenberg Lie algebra instead
of $sl_2.$ The graded algebra $T(\Omega)$ of theta functions
is embedded in the
graded algebra $A(\Omega)$ of auxiliary theta functions in $(Z,W)$ with
$Z,W\in \BC^{(h,g)}$ with respect to $\Omega$ satisfying the following
conditions\,:
\vskip 0.2cm
\ \ $1^0.$ A realization $\left\{\,{\Cal E}_{kl},\,{\Cal D}_{ma},\,\Delta_{nb}
\,\vert\ 1\leq k,l,m,n\leq h,\ 1\leq a,b\leq g\,\right\}$ (cf.\,see section
2 for detail) of the Heisenberg Lie algebra acts on $A(\Omega)$ as
derivations.
\vskip 0.2cm
\ \ $2^0.$ $T(\Omega)$ is the subalgebra consisting of all the elements
$\varphi\in A(\Omega)$ such that $\Delta_{nb}\varphi=0$ for all
$1\leq n\leq h,\ 1\leq b\leq g.$
\vskip 0.2cm
\ \ $3^0.$ The set
$$\left\{\,\Dt\,\bigg|\ \J,\ A\in \LM,\ \M\in \BM\,\right\}$$
forms a canonical basis of $A(\Omega)$.
\vskip 0.2cm
\ \ $4^0.$ The mapping
$$\Dt\lmt \Jt$$
(\,$\J,\ A\in \LM$ and $\M\in \BM$\,) induces an algebra isomorphism of
$A(\Omega)$ onto the algebra of differential polynomials of theta functions.
\vskip 0.3cm
\ \ \
{\smc Notations:}
\ \ We denote by $\BZ,\,\BR$ and $\BC$ the ring of integers, the field of real
numbers, and the field of complex numbers respectively.
The symbol ``:='' means that the expression on the
right is the definition of that on the left. We denote by $\BZ^+$
the set of all positive integers. $F^{(k,l)}$ denotes the set of
all $k\times l$ matrices with entries in a commutative ring $F$.
For any $M\in F^{(k,l)},\ ^t\!M$ denotes the transpose matrix of $M$.
For $A\in F^{(k,k)},\ \sigma(A)$ denotes the trace of $A$. For
$A\in F^{(k,l)}$ and $B\in F^{(k,k)},$ we set $B[A]={^t\!A}BA$.
For a positive symmetric, even integral matrix $\M$ of degree $h,\ \LM$
denotes a complete system of representatives of the cosets
$\M^{-1}\BZ^{(h,g)}/\BZ^{(h,g)}.$
$$\align
\BZ^{(h,g)}_{\geq 0}=&\left\{\,J=(J_{ka})\in \BZ^{(h,g)}\,
\vert\ J_{ka}\geq 0
\ \text{ for\ all}\ k,a\,\right\},\\
\vert J\vert=&\sum_{k,a}\,J_{k,a},\\
J\pm \epsilon_{ka}=&(J_{11},\cdots,J_{ka}\pm 1,\cdots,J_{hg}),\\
J!=&J_{11}!\cdots J_{ka}!\cdots J_{hg}!.
\endalign$$
For $J=(J_{ka})\in \BZ^{(h,g)}_{\geq 0},\,Z=(Z_{ka})$ and $W=(W_{ka}),$
we set
$$Z^J=Z_{11}^{J_{11}}\cdots Z_{hg}^{J_{hg}},\ \ \
W^J=W_{11}^{J_{11}}\cdots W_{hg}^{J_{hg}}.$$
For $J=(J_{ka})\in \BZ^{(h,g)}_{\geq 0},$ we put
$$\left({{\partial}\over {\partial W}}\right)^J=\left(
{{\partial}\over {\partial W_{11}}}\right)^{J_{11}}\cdots \left(
{{\partial}\over {\partial W_{hg}}}\right)^{J_{hg}}.$$
\head 1 \ Auxiliary\ theta\ functions 
\endhead                              
\ \ \ We fix an element $\Omega$ of $H_g$ once and for all. Let $\M$ be
a positive symmetric, even integral matrix of degree $h$. An
{\it auxiliary\ theta\ function\ of\ level\ $\M$} with respect to $\Omega$
means a function $\varphi(Z,W)$ in complex variables $(Z,W)\in \BC^{(h,g)}
\times \BC^{(h,g)}$ such that
\vskip 0.15cm
\ \ (a) $\varphi(Z,W)$ is a polynomial in complex variables $Z=(Z_{ka})$
whose coefficients are entire functions in $W=(W_{ka}),$ and
\vskip 0.15cm
\ \ (b) for all $(\xi,\eta)\in \BZ^{(h,g)}\times \BZ^{(h,g)}$ and
$(Z,W)\in \BC^{(h,g)}\times \BC^{(h,g)},$
$$\varphi(Z+\xi,W+\xi\Omega+\eta)=\emw\,\varphi(Z,W)$$
holds.
\vskip 0.2cm
\ \ \ Let $A_{\M}(\Omega)$ be the vector space of auxiliary theta functions
of level $\M$ with respect to $\Omega.$ We let
$$A(\Omega):=\sum_{\M}\,A_{\M}(\Omega)$$
the graded algebra of auxiliary theta functions, where $\M=(\M_{kl})\
(\,1\leq k,l\leq h\,)$ runs over the set $\BM$ of all positive symmetric,
even integral $h\times h$ matrices such that $\M_{kl}\not= 0$ for all
$k,l.$ We note that $A(\Omega)$ contains the graded algebra $T(\Omega)$
as the subalgebra of polynomials of degree zero in $Z$.
\vskip 0.15cm
\ \ \ We define the auxiliary theta series
$$\align
\ \ \ \ \ &\mjz\\
(1.1)\ \ \ \ \ \ \ \ \ \
=&(2\pi i)^{\vert J\vert}\sum_{\N}\,\prod_{k=1}^h\prod_{a=1}^g\,
\left(\,\sum_{l=1}^h\,\M_{kl}(Z+N+A)_{la}\,\right)^{J_{ka}} \hskip 4cm\\
&\times \text {exp}\left\{ \pi i\s(\M((N+A)\Omega\,^t(N+A)
+2W\,^t(N+A)))\right\},
\endalign$$
where $J=(J_{ka})\in \BZ^{(h,g)}_{\geq 0},\ \M=(\M_{kl})\in \BM,\
A\in \LM$ and $(Z+N+A)_{la}=Z_{la}+N_{la}+A_{la}.$
\vskip 0.2cm
\indent
{\smc Lemma\ 1.1.} For each $\J$ and $\M\in \BM,$ we have
$$\align
\ \ \ &\ {\tilde {\vartheta}}^{(\M)}_J\left[ \matrix A\\ 0\endmatrix\right]
(\Omega\vert Z+\xi,W+\xi\Omega+\eta)\\
\ \ = &\text {exp}\left\{ -\pi i\s(\M(\xi\Omega\,^t\xi+2W\,^t\xi))\right\}\,
\mjz,
\endalign$$
where $A\in \LM$ and $(\xi,\eta)\in \BZ^{(h,g)}\times \BZ^{(h,g)}.$
In particular, $\mjz\in A_{\M}(\Omega).$
\vskip 0.15cm
\indent
{\it Proof.} We observe that
$$\align
\ \ \ \ & (N+A)\Omega\,^t(N+A)+2(W+\xi\Omega+\eta)\,^t(N+A)\\
\ \ \ = & (N+\xi+A)\Omega\,^t(N+\xi+A)+2W\,^t(N+\xi+A)-(N+A)\Omega\,^t\xi\\
\ \ \ \ & + \xi\Omega\,^t(N+A)+2\eta\,^t(N+A)-(\xi\Omega\,^t\xi+2W\,^t\xi).
\endalign$$
In addition, it is easy to see that
$$\s(\M(N+A)\Omega\,^t\xi)=\s(\M\xi\Omega\,^t(N+A))$$
and
$$\s(\M\eta\,^t(N+A))=\s(\M\eta\,^tN)+\s(\M A\,^t\eta)\in \BZ\ \
(\text {because}\ A\in \LM).$$
Therefore the proof follows immediately from these facts.
\hfill $\square$
\vskip 0.3cm
\ \ \ For any $k,a\in \BZ^+$ with $1\leq k\leq h,\ 1\leq a\leq g$ and
$\M=(\M_{kl})\in \BM,$ we put
$$\partial(\M,Z,W)_{ka}:=2\pi i\sum_{l=1}^h\,\M_{kl}Z_{la}\,+\,
{{\partial}\over {\partial W_{ka}}}.\tag 1.2$$
For each $J=(J_{ka})\in \BZ^{(h,g)}_{\geq 0},$ we put
$$\partial(\M,Z,W)^J:=\partial(\M,Z,W)_{11}^{J_{11}}\cdots
\partial(\M,Z,W)_{ka}^{J_{ka}}\cdots
\partial(\M,Z,W)_{hg}^{J_{hg}}.\tag 1.3$$
\indent
Then we obtain the following.
\vskip 0.2cm
\indent
{\smc Lemma\ 1.2.} For each $\J,\ \M\in \BM$ and $A\in \LM,$ we have
$$\mjz=\partial(\M,Z,W)^J\ta.\tag 1.4$$
\indent
{\it Proof.} It is easy to compute that if $\M=(\M_{kl})\in \BM,$
$$\align
\ \ \ & \partial(\M,Z,W)_{ka}\ta\\
\ \ = &\,2\pi i \sum_{\N}\,\left(\,\sum_{l=1}^h\,\M_{kl}(Z+N+A)_{la}\right)\\
\ \ \ &\times \,\text {exp}\left\{2\pi i \s(\M((N+A)\Omega\,^t(N+A)+
2W\,^t(N+A)))\right\}.\endalign $$
The proof follows immediately from the fact that if $J=(J_{ka})\in
\BZ^{(h,g)}_{\geq 0}$ and $\M=(\M_{kl})\in \BM,$
$$\align
\ \ \ & \partial(\M,Z,W)_{ka}^{J_{ka}}\ta\\
\ \ = & (2\pi i)^{J_{ka}}\sum_{\N}\,\left(\,\sum_{l=1}^h\,\M_{kl}
(Z+N+A)_{la}\right)^{J_{ka}}\\
\ \ \ &\times\,\text {exp}\left\{ \pi i\s(\M((N+A)\Omega\,^t(N+A)+
2W\,^t(N+A)))\right\}.\endalign$$
\hfill $\square$
\par\smallpagebreak
\indent
{\smc Theorem\ 1.} For a fixed $\M\in \BM,$ the set
$$\left\{\,\mjz\,\bigg|\ \J,\ A\in \LM\,\right\}$$
is a basis of the vector space $A_{\M}(\Omega)$ of auxiliary theta functions
of level $\M$ with respect to $\Omega.$
\vskip 0.1cm
\indent
{\it Proof.} According to Lemma 1.1, the functions $\mjz\ (\,\J\ \text
{and}\ A\in \LM\,)$ are contained in $A_{\M}(\Omega)$ and it is obvious
that they are linearly independent. We give an ordering $\prec$ on
$\BZ^{(h,g)}_{\geq 0}$ as follows. For $J,K\in \BZ^{(h,g)}_{\geq},$ we
write $J\prec K$ if $\vert J\vert\leq \vert K\vert.$ We say that $Z^K$
has higher degree in $Z$ than $Z^J$ if $J\prec K.$ Now we let
$\varphi(Z,W)=\sum_J\,Z^J f_J(W)$ be an element of $A_{\M}(\Omega)$ and
let $Z^Kf_K(W)$ be one of terms with highest degree $K$ in $Z$.
Since $\varphi\in A_{\M}(\Omega),$ we obtain for each $\xi,\eta\in
\BZ^{(h,g)}$
$$\align
\ \ \ \ \ & \sum_J\,(Z+\xi)^J\,f_J(W+\xi\Omega+\eta)\\
\ \ \ \ = & \emw\,\sum_J\,Z^Jf_J(W).
\endalign$$
Comparing the coefficients of $Z^K,$ we get
$$f_K(W+\xi\Omega+\eta)=\emw\,f_K(W)$$
for each $(\xi,\eta)\in \BZ^{(h,g)}\times \BZ^{(h,g)}.$ Thus $f_K\in
T_{\M}(\Omega)$ and so we obtain
$$f_K(W)=\sum_{\alpha=1}^{(\text {det}\,\M)^g}\,c_{\alpha}\,
\vartheta^{(\M)}\left[ \matrix A_{\alpha}\\ 0\endmatrix\right]
(\Omega\vert W),$$
where $c_{\alpha}\in \BC$ and $A_{\alpha}\in \LM.$ Therefore for suitable
constants $d_{\alpha}\in \BC\ (\,\alpha=1,\cdots,
(\,\text {det}\,\M\,)^g\,)$, the function
$$\varphi(Z,W)-\sum_{\alpha=1}^{(\text {det}\,\M)^g}\,d_{\alpha}\,
{\tilde {\vartheta}}^{(\M)}_K\left[ \matrix A_{\alpha}\\ 0\endmatrix\right]
(\Omega\vert Z,W)$$
is an element of $A_{\M}(\Omega)$ without $Z^K$-term and all the new terms
are of lower degree than $K$ in $Z.$ Continuing this process successively,
we can express $\varphi(Z,W)$ as a linear combination of auxiliary theta
functions
${\tilde {\vartheta}}^{(\M)}_J\left[\matrix A_{\alpha}\\ 0\endmatrix\right]
(\Omega\vert Z,W)\ (\,\J,\ A_{\alpha}\in \LM\,).$
\hfill $\square$
\vskip 0.75cm
\head 2  \ A\ realization\ of Heisenberg\ Lie\ algebra 
\endhead            
\vskip 0.3cm
\ \ \ For each $\M\in \BM,$ we let
$$\s^{(\M)}:A(\Omega)\lrt A_{\M}(\Omega)$$
be the projection operator of $A(\Omega)$ onto $A_{\M}(\Omega).$ We define
the differential operators
$$\align
{\Cal E}_{kl}:=&\sum_{\M\in \BM}\,\M_{kl}\,\s^{(\M)},\ \ \ \M=(\M_{kl}),\\
{\Cal D}_{ma}:=&\sum_{\M\in \BM}\,{1\over {2\pi i}}\,
{{\partial}\over {\partial Z_{ma}}}\circ \s^{(\M)}\,,\\
\Delta_{nb}:=&\sum_{\M\in \BM}\,\partial(\M,Z,W)_{nb}\circ \s^{(\M)},
\endalign$$
where $1\leq k,l,m,n\leq h$ and $1\leq a,b\leq g.$
\vskip 0.15cm
\ \ \ For $J=(J_{ka})\in \bhg,$ we put
$${\Cal D}^J:={\Cal D}_{11}^{J_{11}}\cdots {\Cal D}_{hg}^{J_{hg}},\ \ \
\Delta^J:=\Delta_{11}^{J_{11}}\cdots \Delta_{hg}^{J_{hg}}.$$

\indent
{\smc Proposition\ 2.1.} Let $\M=(\M_{kl})\in \BM$ and $J=(J_{ka})\in
\bhg.$ Then we have
$${\Cal D}_{ma}\mjz=\sum_{l=1}^h\,\M_{ml}J_{la}\,
{\tilde {\vartheta}}^{(\M)}_{J-\epsilon_{la}}\left[\matrix A\\ 0\endmatrix
\right](\Omega\vert Z,W),\tag 2.1$$
$$\Delta_{nb}\mjz={\tilde {\vartheta}}^{(\M)}_{J+\epsilon_{nb}}
\left[\matrix A\\ 0\endmatrix\right](\Omega\vert Z,W),\tag 2.2$$
$$\mjz=\Delta^J\ta,\tag 2.3$$
where $A\in \LM,\ 1\leq m,n\leq h$ and $1\leq a,b\leq g.$
\vskip 0.3cm
\indent
{\it Proof.} (2.1) follows from a direct application of ${\Cal D}_{ma}$ to
(1.1). (2.2) and (2.3) follows immediately from Lemma 1.2, (1.4).
\hfill $\square$
\vskip 0.3cm
\indent {\smc Proposition\ 2.2.} ${\Cal E}_{kl},\ {\Cal D}_{ma},\
\Delta_{nb}\ (\,1\leq k,l,m,n\leq h,\ 1\leq a,b\leq g\,)$ are derivations of
$A(\Omega)$ such that
$$\align
[{\Cal E}_{kl},{\Cal D}_{ma}]&=[{\Cal E}_{kl},\Delta_{nb}]=
[{\Cal D}_{ma},{\Cal D}_{nb}]=[\Delta_{ma},\Delta_{nb}]=0,\\
[{\Cal D}_{ma},\Delta_{nb}]&=\delta_{ab}\,{\Cal E}_{mn}.
\endalign$$

\indent {\it Proof.} According to Proposition 2.1, ${\Cal E}_{kl},\,
{\Cal D}_{ma},\,\Delta_{nb}\,(\,1\leq k,l,m,n\leq h,\,1\leq a,b\leq g\,)$ map
$A(\Omega)$ into itself. Since $A(\Omega)=\sum_{\M\in \BM}\,
A_{\M}(\Omega)$ is a graded algebra, ${\Cal E}_{kl},\,{\Cal D}_{mb},\,
\Delta_{nb}\,(\,1\leq k,l,m,n\leq h,\ 1\leq a,b\leq g\,)$ are derivations
of $A(\Omega)$. An easy calculation yields the above commutation relations.
\hfill $\square$
\vskip 0.3cm
\indent
{\sl Remark\ 2.3.} According to Proposition 2.2,
$\left\{\,{\Cal E}_{kl},\,{\Cal D}_{ma},\,
\Delta_{nb}\,\vert\ 1\leq k,l,m,n\leq h,\ 1\leq a,b\leq g\,\right\}$ is a
realization of the Heisenberg Lie algebra acting on $A(\Omega)$ as
derivations\,(\,cf.\,[Y1],\,[Y2]\,).
\vskip 0.3cm
\indent
{\smc Proposition\ 2.4.} The graded algebra $T(\Omega)$ of theta functions
with respect to $\Omega$ is the subalgebra
of $A(\Omega)$ consisting of $\varphi$ in
$A(\Omega)$ such that ${\Cal D}_{ma}\varphi=0\ (\,1\leq m\leq h,\ 1\leq
a\leq g\,).$
\vskip 0.2cm
\indent
{\it Proof.} If $\varphi\in T(\Omega),\ \varphi$ does not contain any
variables $Z_{ma}$ and ${\Cal D}_{ma}\varphi=0$ for any $m,a\in \BZ^+$
with $1\leq m\leq h,\ 1\leq a\leq g.$\par
Conversely, we assume that
$${\Cal D}_{ma}\left(\,\sum_{J,\alpha,\M}\,c_{J,\alpha,\M}\,
{\tilde {\vartheta}}^{(\M)}_J\left[\matrix A_{\alpha,\M}\\ 0\endmatrix
\right](\Omega\vert Z,W)\right)=0,$$
where $A_{\alpha,\M}\in \LM,\ 1\leq m\leq h$ and $1\leq a\leq g.$ Then
by (2.1), we get
$$\sum_{J,\alpha,\M}\sum_{l=1}^h\,\M_{ml}J_{la}\,c_{J,\alpha,\M}\,
{\tilde {\vartheta}}^{(\M)}_{J-\epsilon_{la}}\left[\matrix
A_{\alpha,\M}\\ 0\endmatrix\right](\Omega\vert Z,W)=0.\tag *$$
Since $\M=(\M_{kl})\in \BM,\ \M_{kl}\not=0$ for all $k,l$ with
$1\leq k,l\leq h.$ Therefore if $J\not=(0,\cdots,0),$ we get
$c_{J,\alpha,\M}=0$ from the condition (*). Hence this completes the proof.
\hfill $\square$
\vskip 0.5cm
\indent
{\smc Theorem\ 2.} $A(\Omega)$ has the direct sum decomposition
$$A(\Omega)=\sum_{\J}\,\Delta^J T(\Omega)=\sum_{\J}\sum_{\M\in \BM}\,
\Delta^J T_{\M}(\Omega) \tag 2.4$$
such that $\Delta^J$ induces a vector space isomorphism of $T_{\M}(\Omega)$
onto $\Delta^J T_{\M}(\Omega).$
\vskip 0.3cm
\indent
{\it Proof.} The proof follows from (2.3) and the fact that
$$\left\{\,\mjz\,\bigg|\ \J,\ A\in \LM,\ \M\in \BM\,\right\},$$
$$\left\{\ \mjz\,\bigg|\ A\in \LM\,\right\}$$
and
$$\left\{\, \ta\,\bigg|\ A\in \LM\,\right\}$$
are the bases of $A(\Omega),\,\Delta^J T_{\M}(\Omega)$ and $T_{\M}(\Omega)$
respectively.
\hfill $\square$
\vskip 0.5cm
\indent
{\sl Remark\ 2.5.}\ We may express the inverse mapping of $\Delta^J:
T_{\M}(\Omega)\lrt \Delta^J T_{\M}(\Omega)$ in terms of ${\Cal D}_{ma},\,
\Delta_{nb}\ (\,1\leq m,n\leq h,\ 1\leq a,b\leq g\,).$ The expression is very
complicated and so we omit it.
\vskip 0.75cm
\head 3 \  Decomposition\ theorem\ on\ differential\ polynomials\\ of\
theta\ functions\ of\ high\ level 
\endhead   
\vskip 0.5cm
\ \ \ In this section, we prove the algebra isomorphism theorem.
\vskip 0.35cm
\indent
{\smc Theorem\ 3.} The replacement
$$\Delta^J \varphi(W)\lrt \left({{\partial}\over {\partial W}}\right)^J
\varphi(W)\ \ \ (\,\J,\,\varphi\in T(\Omega)\,)$$
induces a $T(\Omega)$-algebra isomorphism of $A(\Omega)$ onto the algebra
$$\CJ,\ \ \ A_{\alpha,\M}\in \LM,\ \M\in \BM$$
of differential polynomials of theta functions, namely
\vskip 0.3cm
\indent
(1)$\ \ \GJ=0$\par
\indent
$\ \ \ $ if and only if $\GW=0,$
\vskip 0.25cm
\indent
(2)$\ \ \GJ=\GW$\par
\indent
$\ \ \ \text {if\ and\ only\ if}\ \GW\in T(\Omega).$
\vskip 0.25cm
\indent
{\it Proof.} It is enough to assume that $\GJ$ belongs to $A_{\M}(\Omega)$
for some $\M\in \BM.$ Suppose $\GJ=0.$ By putting $Z=0,$ we obtain
$$\GW=0.$$
Conversely, we suppose that $\GW=0.$ According to Theorem 2, we may write
$$\GJ=\sum_K\,\Delta^K \phi_K(W),\tag 3.1$$
where $\phi_K\in T_{\M}(\Omega).$ Then we have
$$\align
\sum_K\,\left({{\partial}\over {\partial W}}\right)^K\,\phi_K(W)
=&\GJ\bigg|_{Z=0}\\
=&\GW=0.\endalign$$
Therefore it suffices to show $\phi_K(W)=0$ under the condition
$$\sum_K\,\left({{\partial}\over {\partial W}}\right)^K\,\phi_K(W)=0
\ \ \text {and}\ \ \phi_K\in T_{\M}(\Omega).$$
For each $\xi\in \BZ^{(h,g)},$ we get
$$\phi_K(W+\xi\Omega)=\emw\,\phi_K(W)$$
and
$$\align
\sum_K\,\left({{\partial}\over {\partial W}}\right)^K\,\phi_K(W+\xi\Omega)
=&\sum_K\,\pw^K\left[\emw\,\phi_K(W)\right] \\
=&\emw\\
 &\times\,\sum_K\sum_P\,\left(\matrix K\\ P\endmatrix\right)\,
 (-2\pi i\,\M\xi)^P\,\pw^{K-P}\phi_K(W).
\endalign$$
Here if $K=(K_{ma})$ and $P=(P_{ma})$ in $\bhg,$ we put
$$\left(\matrix K\\ P\endmatrix\right):=\left(\matrix K_{11}\\P_{11}
\endmatrix\right)\cdots \left(\matrix K_{ma}\\P_{ma}\endmatrix\right)\cdots
\left(\matrix K_{hg}\\ P_{hg}\endmatrix\right)$$
and if $\M=(\M_{kl})$ and $\xi=(\xi_{ma})\in \BZ^{(h,g)},$ we put
$$(-2\pi i\,\M\xi\,)^P:=\left(-2\pi i\sum_{l=1}^h\,\M_{1l}\xi_{l1}
\right)^{P_{11}}\cdots \left(-2\pi i\sum_{l=1}^h\,\M_{hl}\xi_{lg}
\right)^{P_{hg}}.$$
Thus we have
$$\sum_K\sum_P\,\left(\matrix K\\ P\endmatrix\right)\,(-2\pi i\,\M\xi\,)^P\,
\pw^{K-P}\phi_K(W)=0 \tag 3.2$$
for all $\xi\in \BZ^{(h,g)}.$ Let $K_0$ be one of maximal $K$ in the above
sum (3.2). Then the coefficients of $\xi^{K_0}$ in the polynomial relation
(3.2) in $\xi$ is given by $C(\M)\,\phi_{K_0}(W)$ with nonzero constant
$C(\M)\not=0.$ Thus we get $\phi_{K_0}(W)=0.$ Continuing this process
successively, we have $\phi_K=0$ for all $K$ appearing in the sum (3.1).
Hence from (3.1), we have
$$\GJ=0.$$
We assume that
$$\GJ=\GW.$$
Then we have, for any $(\xi,\eta)\in \BZ^{(h,g)}\times \BZ^{(h,g)},$
$$\align
 &\GW\bigg|_{W\mapsto W+\xi\Omega+\eta}\\
=&\GJ\bigg|_{(Z,W)\mapsto (Z+\xi,W+\xi\Omega+\eta)}\\
=&\emw\,\GJ\\
=&\emw\,\GW.
\endalign$$
Therefore we obtain
$$\GW\in T_{\M}(\Omega).$$
Conversely, we assume that
$$\GW\in T_{\M}(\Omega).$$
Applying (1) to
$$\align
\ \ & F\left(\cdots,\Delta^J
\vartheta^{(\M)}\left[\matrix A_{\alpha,\M}\\ 0\endmatrix\right](\Omega\vert
W),\cdots\right)\\
=&\GJ-\GW, \endalign$$
we obtain
$$F\left(\cdots,\Delta^J \vartheta^{(\M)}\left[\matrix A_{\alpha,\M}\\
0\endmatrix\right](\Omega\vert W),\cdots\right)=0.$$
Hence we get
$$\GJ=\GW.$$
\hfill $\square$
\vskip 0.2cm
\indent
Combining Theorem 2 and Theorem 3, we obtain the decomposition theorem.
\vskip 0.5cm
\indent
{\smc Theorem\ 4.} The algebra $\CJ$ of differential polynomials of theta
functions has a canonical linear basis
$$\left\{\,\pw^J\vartheta^{(\M)}\left[\matrix A_{\alpha,\M}\\
0\endmatrix\right](\Omega\vert W)\,\bigg|\ \J,\ A_{\alpha,\M}\in \LM,
\ \M\in \BM\,\right\},\tag 3.3$$
namely, differential polynomials of theta functions are uniquely expressed
as linear combinations of (3.3) with constant coefficients depending only
on $\Omega.$
\vskip 0.5cm
\indent
{\sl Remark\ 3.1.} In [M3], Morikawa proved the decomposition theorem on
differential polynomials of theta functions in the case that $h=1.$ He
investigated the graded algebras of theta functions and of auxiliary
theta functions\,:
$$\Theta_0(\Omega)=\sum_{n=1}^{\infty}\,\Theta_0^{(n)}(\Omega),\ \ \
\ \Theta(\Omega)=\sum_{n=1}^{\infty}\,\Theta^{(n)}(\Omega),$$
where $\Theta_0^{(n)}(\Omega)$\,(\,resp.\,
$\Theta^{(n)}(\Omega)\,)$ denotes the
vector space of theta functions\,(\,resp.\,auxiliary theta functions\,) of
level $n$ with respect to $\Omega$. In this paper, when $h=1,$ we
investigated the following graded algebras
$$T(\Omega)=\sum_{n=1}^{\infty}\,\Theta_0^{(2n)}(\Omega),\ \ \ \
A(\Omega)=\sum_{n=1}^{\infty}\,\Theta^{(2n)}(\Omega)$$ of theta
functions and of auxiliary theta functions of {\it even\ level}
with respect to $\Omega.$

\par\bigpagebreak
Department of Mathematics\par Inha University\par Incheon
402-751\par Republic of Korea \vskip 1cm
email\,:\,jhyang\@inha.ac.kr

\end
\vskip 1cm \Refs \widestnumber\key{\bf M3}






\ref \key {\bf M1}  \by H. Morikawa \paper Some analytic and
geometric applications of the invariant theoretic method \jour
Nagoya Math. J. \vol 80 \yr 1980 \pages 1-47 \endref

\ref\key {\bf M2} \bysame \book On Possion brackets of
semi-invariants, Manifolds and Lie groups \publ Progress in Math.
Birkh{\"a}user  \yr 1981 \pages 267-281 \endref

\ref\key{\bf M3}\bysame \paper A decomposition theorem on
differential polynomials of theta functions \jour Nagoya Math. J.
\vol 96 \yr 1984 \pages 113-126 \endref

\ref\key{\bf Y1}  \by J.-H. Yang \paper Harmonic analysis on the
quotient spaces of Heisenberg groups \jour Nagoya Math. J. \vol
123 \yr 1991 \pages 103-117
\endref

\ref\key {\bf Y2}\bysame \paper Harmonic analysis on the quotient
spaces of Heisenberg groups, II \jour J. Number Theory \vol 48 \yr
1994  \endref
\endRefs
